\documentclass[12pt]{article}
\usepackage{fullpage}
\usepackage{amssymb}
\usepackage{amsthm}
\usepackage{graphicx}
\usepackage[all,2cell]{xy}
\newtheorem{theorem}{Theorem}[section]
\newtheorem{corollary}[theorem]{Corollary}

\newtheorem{lemma}[theorem]{Lemma}

\begin{document}

\title{Projective normality of finite group quotients and EGZ theorem}
\author{S.S.Kannan, S.K.Pattanayak \\  
\\ Chennai Mathematical Institute, Plot H1, SIPCOT IT Park,\\ Padur 
Post Office, Siruseri, Tamilnadu - 603103, India.\\
kannan@cmi.ac.in, santosh@cmi.ac.in, pranab@cmi.ac.in} 

\maketitle
\date{}

\begin{abstract} In this note, we prove that for any finite dimensional vector 
space $V$ over $\mathbb {C}$, and for a finite cyclic group $G$, the projective variety $\mathbb P(V)/G$ is projectively normal with respect to the descent of $\mathcal O(1)^{\otimes |G|}$ by a method using toric variety, and deduce the EGZ theorem as a consequence. 
 
\end{abstract}
\hspace*{2.5cm}Keywords: GIT quotient, line bundle, normality of a semigroup. 
\begin{center}
{\bf Introduction}
\end{center}

Let $V$ be a finite dimensional representation of a cyclic group $G$ over the field of complex numbers $\mathbb {C}$. Let $\mathcal {L}$ denote the descent of $\mathcal {O}(1)^{\otimes |G|}$ to the GIT quotient $\mathbb {P}(V)/G$. In [4], it is shown that $(\mathbb {P}(V)/G, \mathcal {L})$ is projectively normal. Proof of this uses the well known arithmetic result due to Erd\"{o}s-Ginzburg-Ziv (see [2]).

In this note, we prove the projective normality of $(\mathbb {P}(V)/G, \mathcal {L})$ by a method using toric variety, and deduce the EGZ theorem (see [2]) as a consequence. 

\section{Erd\"{o}s-Ginzburg-Ziv theorem:}
We first prove a lemma on normality of a semigroup related to a finite cyclic group.

\begin{lemma}
Let $M$ be the sub-semigroup of $\mathbb {Z}^{n}$ generated by the finite set $S= \{(m_0,m_1, \cdots ,m_{n-1}) \in (\mathbb {Z}_{\geq 0})^{n}: \sum_{i=0}^{n-1}m_i=n \,\, and \,\, \sum_{i=0}^{n-1}im_i \equiv 0 \,\, mod \,\, n\}$ and let $N$ be the subgroup of $\mathbb {Z}^{n}$ generated by $M$. Then $M = \{x \in N: qx \in M \,\, for \,\, some \,\, q \in \mathbb {N}\}$.
\end{lemma}

\begin{proof}
Consider the homomorphism:
            
$\Phi: \mathbb{Z}^{n}\longrightarrow \frac {\mathbb{Z}}{n\mathbb{Z}}\oplus \frac {\mathbb{Z}}{n\mathbb{Z}}$ of abelian groups given by:

 $\Phi(x_0,x_1,\cdots ,x_{n-1})= (\sum_{i=0}^{n-1}x_i+n \mathbb{Z},  \sum_{i=0}^{n-1}ix_i+n \mathbb{Z})$. 

Clearly, $\Phi$ is surjective and $N \subset Ker(\Phi)$. 
So $\frac {\mathbb{Z}^{n}} {Ker (\Phi)} \cong \frac {\mathbb{Z}}{n\mathbb{Z}}
\oplus \frac {\mathbb{Z}}{n\mathbb{Z}} \hspace {2cm} \longrightarrow (1)$.

Now, we show that $N = Ker(\Phi)$.

Let $\{e_i, i = 0,1,2, \cdots ,n-1\}$ be the standard basis of $\mathbb {Z}^{n}$. Then the subgroup of $N$ generated by $\{u_{n-2}=ne_{n-2}, u_{n-1}=ne_{n-1}, u_r= e_r+(r+1)e_{n-2}+(n-r-2)e_{n-1}, r= 0,1,2, \cdots ,n-3\}$ is of index $n^2$ in $\mathbb {Z}^{n}$.

On the otherhand, $N \subset Ker(\Phi)$. Hence, $N = Ker(\Phi)$, by (1). 

Now each $(m_0,m_1, \cdots ,m_{n-1}) \in M$ can be written as a $\mathbb {Z}$-linear combination of $u_i$'s: $(m_0,m_1, \cdots ,m_{n-1})= \sum_{i=0}^{n-1}d_iu_i$, where $d_i= m_i \in \mathbb {Z}_{\geq 0}, i= 0,1,2, \cdots ,n-3$, $d_{n-2}= \frac {m_{n-2}- \sum_{i=0}^{n-3}(i+1)m_i} {n}$ and $d_{n-1}=  \frac {m_{n-1}- \sum_{i=0}^{n-3}(n-i-2)m_i}{n}$. Notice that $d_{n-2}, d_{n-1} \in \mathbb {Z}$ by the conditions on $N$. 

Let $x \in N$ be such that $qx \in M$, for some $q \in N$.

Then $q(\sum_{i=0}^{n-1}a_iu_i)= \sum_{i=0}^{n-1}b_iu_i+ \sum_{j=1}^{m}c_jv_j$, where $\{v_j:j= 1,2,\cdots ,m\}= S \setminus  \{u_i:i= 0,1,2,\cdots ,n-1\}, a_i \in \mathbb {Z}, b_i, c_j \in \mathbb {Z}_{\geq 0} \,\, \forall \,\, i, j$.

Again, we can write $v_k=  \sum_{i=0}^{n-1}d_{i,k}u_i, d_{i,k}\in \mathbb {Z}\,\, \forall \,\, i, d_{i,k}\in \mathbb {Z}_{\geq 0} \,\, \forall \,\, i= 0,1,2, \cdots ,n-3$ and $\exists \,\, 0 \leq i \leq (n-3)$ such that $d_{i,k} > 0$.

If $x \notin M$, we may assume that $a_i \leq 0 \,\,\forall \,\,i$.

If one of the $b_i$'s or $c_j$'s is nonzero, then there is an $i$ for which $a_i > 0$, contradiction to the assumption that  $a_i \leq 0$. So, $x \in M$.

\end{proof}

We now prove:

\begin{theorem}
Let $G$ be a cyclic group of order $n$ and $V$ be any finite dimensional representation of $G$ over $\mathbb {C}$. Let $\mathcal {L}$ be the descent of $\mathcal {O}(1)^{\otimes n}$. Then $(\mathbb {P}(V)/G, \mathcal {L})$ is projectively normal. 

\end{theorem}

\begin{proof}

Let $R:=\oplus_{d\geq 0}R_d$; $R_d:= (Sym^{dn}V)^G$. 

Let $G=<g>$.  Write $V^*=\oplus_{i=0}^{n-1}V_i$ 
where $V_i:=\{ v\in V^*: g.v=\xi^i.v\}$, $0\leq i\leq n-1$, where $\xi$ is a primitive $n$th root of unity. The $\mathbb {C}$-vector space $R_1$ is generated by elements of the form $X_0.X_1\ldots X_{n-1}$, where $X_i \in Sym^{m_i}(V_i), \sum_{i=0}^{n-1}m_i=n \,\, and \,\, \sum_{i=0}^{n-1}im_i \equiv 0 \,\, mod \,\, n$. 

So, the $\mathbb {C}$-subalgebra of  $\mathbb {C}[V]$ generated by $R_1$ is the algebra corresponding to the semigroup $M$ generated by $\{(m_0,m_1, \cdots ,m_{n-1}) \in (\mathbb {Z}_{\geq 0})^{n}: \sum_{i=0}^{n-1}m_i = n \,\, and \,\, \sum_{i=0}^{n-1}im_i \equiv 0 \,\, mod \,\, n\}$.

By lemma (1.1), $M$ is normal (for the definition of normality of a semigroup, see page 61 of [1]).

Hence, by theorem 4.39 of [1] the $\mathbb {C}$-subalgebra of  $\mathbb {C}[V]$ generated by $R_1$ is normal.

Thus, by Exercise 5.14(a) of [3], the theorem follows. 

\end{proof}

We now deduce EGZ-theorem.

\begin{corollary}
Let $\{a_1,a_2,\cdots ,a_m\}, m \geq 2n-1$ be a sequence of elements of $ \frac {\mathbb {Z}}{n\mathbb {Z}}$. Then there exists a subsequence $\{a_{i_1}, a_{i_2}, \cdots ,a_{i_n}\}$ of length $n$ whose sum is zero.
\end{corollary}

\begin{proof}
Let $G=\frac {\mathbb {Z}}{n\mathbb {Z}} = <g>$ and $V$ be the regular representation of $G$ over $\mathbb {C}$. 

Let $\{X_i: i=0,1,\cdots ,n-1\}$ be a basis of $V^*$ given by:

$\hspace {2cm} g.X_i= \xi^iX_i, \,\, \forall \,\, g \in G$ and $i=0,1,\cdots ,n-1$, where $\xi$ is a primitive $n$th root of unity.

Let $\{a_1, a_2, \cdots ,a_{m}\}, m \geq 2n-1$ be a sequence of elements of  $G$. Consider the subsequence $\{a_1, a_2, \cdots ,a_{2n-1}\}$ of length $2n-1$.

Take $a = -(\sum_{i=1}^{2n-1}a_i)$.

Then $(\prod_{i=1}^{2n-1}X_{a_i}).X_{a}$ is a $G$-invariant monomial of degree $2n$.

By Theorem (1.2), there exists a subsequence  $\{a_{i_1}, a_{i_2}, \cdots ,a_{i_n}\}$ of  $\{a_1, a_2, \cdots ,a_{2n-1}, a\}$ of length $n$ such that $\prod_{j=1}^{n}X_{a_{i_j}}$ is $G$-invariant.

So, $\sum_{j=1}^{n}a_{i_j} = 0$. Hence, the Corollary follows.

\end{proof}


\begin{thebibliography}{22}

\bibitem[1]{r1} W.Bruns, J.Gubeladze, Polytopes, Rings, and K-Theory. Springer Monographs
in Mathematics. Springer. to appear.
\bibitem[2]{r2} P.Erd\"{o}s, A.Ginzburg, A.Ziv, A theorem in additive number theory,
Bull. Res. Council, Israel, 10 F(1961) 41-43. 
\bibitem[3]{r3} R.Hartshorne, Algebraic Geometry, Springer-Verlag, 1977.
\bibitem[4]{r4} S.S.Kannan, S.K.Pattanayak, Pranab Sardar, Projective normality of finite groups quotients. Proc. Amer. Math. Soc. 137 (2009), no. 3, pp. 863-867. 

\bibitem[5]{r5} D.Mumford, J.Fogarty and F.Kirwan, Geometric Invariant theory,
Springer-Verlag, 1994.



\end{thebibliography}
\end{document}